\documentclass[12pt]{amsart}
\usepackage{amsxtra}
\usepackage{url}
\usepackage[utf8]{inputenc}
\usepackage{hyperref}
\newtheorem{definition}{Definition}
\newtheorem{state}{Proposition}

\newtheorem{cor}{Corollary}
\newtheorem{note}{Note}
\newtheorem{lem}{Lemma}
\usepackage{mathtools}

\DeclarePairedDelimiter\floor{\lfloor}{\rfloor}

\newcommand{\fl}[1]{\floor*{\frac{#1}{2}}}
\newcommand{\qnum}[1]{\begin{bmatrix}
#1
\end{bmatrix}_q
}
\newcommand{\qbinom}[2]{\begin{bmatrix}
#1 \\
#2
\end{bmatrix}_q
}
\newcommand{\qbinomm}[2]{\begin{bmatrix}
#1 \\
#2
\end{bmatrix}_{q^2}
}

\title[$q$-$\gamma$-positivity of two $q$-Narayana polynomials.]%
{$q$-$\gamma$-Positivity of two $q$-Narayana polynomials.}

\begin{document}

\author{Lenny Tevlin}
\address[Lenny Tevlin]{No longer affiliated with an academic institution}
\email{ltevlin@gmail.com}

\begin{abstract}
We address the question of  $q$-$\gamma$-positivity of $q$-Narayana polynomials for $q$-super Catalan numbers $A_{n,m}(q)$ and their duals, $B_{n,m}(q)$.

We offer three main results.
First, both the $q$-$m$-Narayana polynomials, whose coefficients refine the $q$-super Catalan numbers, and the dual $q$-Narayana polynomials are $q$-$\gamma$-positive of index one. Second, the dual $q$-Narayana polynomials are also $q$-$\gamma$-positive of index two. Third, both the  $q$-Narayana numbers and their duals are coefficient-wise $q$-unimodal. 
\end{abstract}

\maketitle
\tableofcontents

\section{Introduction}
In 2014 Michelle Wachs gave a talk at the IMA workshop Geometric and Enumerative Combinatorics, where she presented a review of $q$-$\gamma$-positivity~\cite{Wachs2014}. 
 We begin by recalling the definition of $q$-$\gamma$-positive polynomials given by Krattenthaler and Wachs in~\cite{Wachs2014}:

\begin{definition}
\label{def:qgamma}
Let $a\ge1$ and $d\ge0$ be integers. A polynomial $F(q,t)\in\mathbb Z_{\ge0}[q,t]$ of degree $d$ in $t$ is $q$-$\gamma$-positive of index $a$ if
\begin{equation}
\label{eq:qgamma-definition}
F(q,t)=\sum_{j=0}^{\lfloor d/2\rfloor} q^{a\binom{j+1}{2}}\gamma_j(q)\,t^j \prod_{i=j+1}^{d-j}(1+q^{ai}t)
\qquad\text{with}\qquad \gamma_j(q)\in\mathbb Z_{\ge0}[q].
\end{equation}
\end{definition}

We are interested in two families of $q$-factorial ratios:
\begin{align}
\label{def:sC}
& A_{n,m}(q) = \frac{\qnum{2m}! \qnum{2n}!}{\qnum{m}! \qnum{n}! \qnum{m +n}!} =  \frac{( q ; q )_{2 m} ( q ; q )_{2 n}}{ ( q ; q )_m ( q ; q )_n ( q ; q )_{n + m}}\qquad  \text{and} \\
\label{eq:B-def}
& B_{n,m}(q)=\frac{\qnum{2n}!\qnum{m}!} {\qnum{n}!\qnum{2m}!\qnum{n-m}!} = \frac{( q ; q )_{2 n} ( q ; q )_{m}}{ ( q ; q )_n ( q ; q )_{2 m}  ( q ; q )_{n - m}},
\qquad 0\le m\le n.
\end{align}
 Warnaar and Zudilin~\cite{Warnaar2011} proved that they are polynomials in $q$ with nonnegative integer coefficients.  We derived a number of combinatorial identities satisfied by~\eqref{def:sC} and studied their $q$-Narayana refinements in~\cite{RodeletTevlin}.

The polynomials $A_{n,m}(q)$ are referred to as the $q$-super Catalan numbers, and we call $B_{n,m}(q)$ the dual $q$-super Catalan numbers because of the relation
\begin{equation}
\label{eq:intro-duality}
A_{n,m}(q)\,B_{n,m}(q)=\qbinom{2n}{n}\qbinom{2n}{n-m},
\end{equation}
which follows by cancelling the factors $\qnum{2m}!$ and $\qnum{m}!$. 
Also, directly from the definitions,
\begin{equation}
\label{eq:AB-ratio}
B_{n,m}(q)=\qbinom{n+m}{2m}\qbinom{2m}{m}^{-1}A_{n,m}(q).
\end{equation}

The dual family includes the $q$-central binomial coefficients, $B_{n,0}(q)=\qbinom{2n}{n}=A_{n,0}(q)$, and $B_{m,m}(q)=1$. Unlike the $q$-super Catalan numbers, the $B_{n,m}(q)$ are not symmetric in $n$ and $m$.

Recall from~\cite{Wachs2014} and~\cite[Eq.~(85)]{RodeletTevlin} that the $q$-super Catalan numbers decompose in a Narayana-like way,
\begin{equation}
\label{eq:A-weighted}
A_{n,m}(q)=\sum_{k=0}^{n-m}q^{k(k+m)}N_{n,k}^{A,(m)}(q),
\end{equation}
where $N_{n,k}^{A,(m)}(q)\in\mathbb Z_{\ge0}[q]$ are the $q$-$m$-Narayana numbers. In~\cite[Eq.~(86)]{RodeletTevlin} they are denoted $N^{(m)}_{n,k}(q)$; here the superscript $A$ distinguishes them from those for the dual family.

The duality relation~\eqref{eq:AB-ratio} is a driving principle throughout the paper.

In Sec.~\ref{sec:qNarB}, we show that it persists at the level of the $q$-Narayana refinements: the dual $q$-super Catalan numbers admit a Narayana-like decomposition, and the two $q$-Narayana families are related by the same factor. This relation converts the positive expansion for the $q$-$m$-Narayana numbers into the corresponding positive expansion for the dual $q$-$m$-Narayana numbers, yielding their integrality and positivity.

Our main results are:
\begin{enumerate}
\item The polynomials $\mathcal{N}^{A,(m)}_n(q,t)=\sum_kq^{\binom{k+1}{2}}N^{A,(m)}_{n,k}(q)t^k$ are $q$-$\gamma$-positive of index $1$
(Proposition~\ref{prop:A-qgamma}).
\item The  polynomials $\mathcal{N}^{B,(m)}_n(q,t)=\sum_kq^{\binom{k+1}{2}}N^{B,(m)}_{n,k}(q)t^k$ are also $q$-$\gamma$-positive of index $1$
(Proposition~\ref{prop:B-qgamma}).
\item
The polynomials $\widetilde{\mathcal{N}}^{B,(m)}_n(q,t)=\sum_kq^{k(k+1)}N^{B,(m)}_{n,k}(q)t^k$ are $q$-$\gamma$-positive of index $2$
(Proposition~\ref{prop:B-index2-qgamma}).
\item 
Every row of $q$-$m$-Narayana numbers and every row of dual $q$-Narayana numbers is coefficient-wise $q$-unimodal (Corollary~\ref{cor:row-q-unimodal}).
\end{enumerate}

Krattenthaler and Wachs~\cite{Wachs2014} proved that $q$-$m$-Narayana numbers are coefficient-wise $q$-unimodal and showed that the polynomials $\sum_k q^{k(k+1)}N^{A,(m)}_{n,k}(q)t^k$, although positive, are not $q$-$\gamma$-positive of index $2$ in general.

\subsection*{Outline of the paper}

After recalling the standard notation for $q$-analogs and $q$-series, we introduce the dual $q$-$m$-Narayana numbers, prove a Narayana-like decomposition of $B_{n,m}(q)$, and prove that the dual $q$-$m$-Narayana numbers are palindromic polynomials in $q$ with nonnegative integer coefficients. 

Next we prove $q$-$\gamma$-positivity of index $1$, first for the $q$-$m$-Narayana polynomials and then for the dual $q$-$m$-Narayana polynomials. We then prove that the dual $q$-Narayana polynomials are $q$-$\gamma$-positive of index $2$.

After that we establish a new $q$-Touchard identity for dual $q$-super Catalans. 

Finally, we deduce the coefficient-wise $q$-unimodality of the rows of both families. We close with the classical case $q=1$.

\subsection*{Notation and $q$-series identities}

We use the following notation:
\begin{equation}
(a;q)_n =
\begin{cases}
\prod_{k=0}^{n-1} (1 - a q^k) & \text{for\ } n \geq 1, \\
1 & \text{for\ } n = 0,
\end{cases}
\end{equation}
and $(a_1,\ldots,a_s;q)_n=(a_1;q)_n\cdots(a_s;q)_n$. The basic hypergeometric series are
\begin{equation}
\label{eq:series}
{}_{r+1}\phi_r\!\left[\begin{matrix} a_1,\ldots,a_{r+1}\\
b_1,\ldots,b_r\end{matrix};q,z\right]
=\sum_{j\ge0}\frac{(a_1,\ldots,a_{r+1};q)_j}{(q,b_1,\ldots,b_r;q)_j}\,z^j .
\end{equation}
Such a series is terminating if one of the $a_i$ equals $q^{-N}$ for an integer $N\ge0$, and, for $z=q$, it is balanced if $q\,a_1\cdots a_{r+1}=b_1\cdots b_r$. Further,
\begin{equation}
\qnum{n}=\frac{1-q^n}{1-q},\qquad \qnum{n}!=\prod_{i=1}^{n}\qnum{i},\qquad
\qbinom{n}{k}  = \frac{(q^{n-k+1};q)_k}{(q;q)_k}
\end{equation}
Furthermore, we will use the following $q$-series identities in proofs:
\begin{align}
\label{eq:Ksplit}
& (a;q)_{2n} = (a;q^2)_n\,(aq;q^2)_n, \\
\label{eq:Ktwo}
& (a^2;q^2)_n = (a;q)_n\,(-a;q)_n, \\
\label{eq:Kfour}
& (\pm q^a;q^b)_n = (\mp 1)^n q^{an+b\binom{n}{2}}\,(\pm q^{-a-b(n-1)};q^b)_n,
\end{align}
where in~\eqref{eq:Kfour} the upper signs and the lower signs are taken together.  We will also use the finite $q$-binomial theorem~\cite[Exercise~1.2(vi)]{GR} in the form
\begin{equation}
\label{eq:q-binomial-theorem}
\prod_{i=k}^{k+n-1}(1+q^i t) =\sum_{\ell=0}^{n}q^{k \ell+\binom{\ell}{2}} \qbinom{n}{\ell}t^\ell
\qquad(n, k \ge0),
\end{equation}
In addition, we write $A_{n,m}=A_{n,m}(1), B_{n,m}=B_{n,m}(1)$ and $N^{A,(m)}_{n,k}=N^{A,(m)}_{n,k}(1), N^{B,(m)}_{n,k}=N^{B,(m)}_{n,k}(1)$, etc. 

\section{Dual $q$-Narayana numbers and their positivity}
\label{sec:qNarB}
We shall introduce the $q$-Narayana refinement of $B_{n,m}(q)$ akin to that of $A_{n,m}(q)$~\cite[Eq.~(86)]{RodeletTevlin}:
\[
A_{n,m}(q)=\sum_{k=0}^{n-m}q^{k(k+m)}N_{n,k}^{A,(m)}(q),
\]
where the $q$-$m$-Narayana numbers are
\begin{equation}
\label{eq:N_def}
N^{A,(m)}_{n,k}(q)  = \qbinom{2m}{m}\qbinom{n}{m}^{-1} \qbinom{n}{k} \qbinom{n}{k+m}, \qquad N^{A,(m)}_{n,k}(q) \in \mathbb Z_{\ge0}[q].
\end{equation}

\begin{definition}
Let $m$, $n$, and $k$ be integers such that $0\le m\leq n$ and $0\le k\leq n-m$.
The dual $q$-$m$-Narayana numbers are
\begin{equation}
\label{eq:B-Nara}
N^{B,(m)}_{n,k}(q)
 = \qbinom{n+m}{m}\qbinom{2m}{m}^{-1} \qbinom{n}{k} \qbinom{n}{k+m},
\end{equation}
or, equivalently,
\begin{equation}
\label{eq:B-Nara-factorial}
N^{B,(m)}_{n,k}(q) =\frac{(q;q)_{n+m}\,(q;q)_m\,(q;q)_n} {(q;q)_{2m}\,(q;q)_{n-m-k}\,(q;q)_k\,(q;q)_{n-k}\,(q;q)_{k+m}}.
\end{equation}
\end{definition}

Comparing~\eqref{eq:N_def} and~\eqref{eq:B-Nara}, and using~\eqref{eq:AB-ratio}, we obtain
\begin{equation}
\label{eq:AB-Nara}
N^{B,(m)}_{n,k}(q)=\qbinom{n+m}{2m}\qbinom{2m}{m}^{-1}N^{A,(m)}_{n,k}(q)=\frac{B_{n,m}(q)}{A_{n,m}(q)}\,N^{A,(m)}_{n,k}(q).
\end{equation}
Thus the dual $q$-$m$-Narayana numbers are dual to the $q$-$m$-Narayana numbers in the same sense that $B_{n,m}(q)$ is dual to $A_{n,m}(q)$.

\begin{state}
\label{prop:B-weighted}
Let $m$ and $n$ be two nonnegative integers such that $m\leq n$. Then the dual $q$-super Catalan numbers decompose in a Narayana-like way as
\begin{equation}
\label{eq:B-weighted}
B_{n,m}(q) = \sum_{k=0}^{n -m} q^{k(k+m)} N_{n,k}^{B,(m)}(q).
\end{equation}
\end{state}

\begin{proof} 
By~\eqref{eq:AB-Nara},~\eqref{eq:A-weighted}, and~\eqref{eq:AB-ratio}, 
\[ 
\sum_{k=0}^{n-m}q^{k(k+m)}N^{B,(m)}_{n,k}(q)=\frac{B_{n,m}(q)}{A_{n,m}(q)}\sum_{k=0}^{n-m}q^{k(k+m)}N^{A,(m)}_{n,k}(q)=\frac{B_{n,m}(q)}{A_{n,m}(q)}A_{n,m}(q)=B_{n,m}(q). 
\] \end{proof}
Both families are symmetric:
\begin{equation}
\label{eq:B-reflection}
N_{n,k}^{A,(m)}(q)=N_{n,n-m-k}^{A,(m)}(q), \qquad N_{n,k}^{B,(m)}(q)=N_{n,n-m-k}^{B,(m)}(q).
\end{equation}
since replacing $k$ by $n-m-k$ interchanges the two $q$-binomial coefficients in~\eqref{eq:N_def} and in~\eqref{eq:B-Nara}.

In~\cite[Eq.~(117)]{RodeletTevlin} we proved the following expansion of the $q$-$m$-Narayana numbers, valid for $0\le k\le\fl{n-m}$:
\begin{equation}
\label{eq:A-positive-expansion}
N^{A,(m)}_{n,k}(q) = \sum_{j=0}^{\min(k,n-m-k)} q^{j(j+m)} \qbinom{n-m}{2j} A_{m,j}(q) \qbinom{n - 2j-m}{k-j},
\end{equation}
which implies their integrality and positivity.

The integrality and positivity of the dual $q$-$m$-Narayana numbers do not follow directly from their definition. Our proof is based on the following analog of this identity.

\begin{state}
\label{prop:B-positive-expansion}
Let $n$, $m$, and $k$ be nonnegative integers such that $m\leq n$ and $k\leq n-m$. Then
\begin{equation}
\label{eq:B-positive-expansion}
N^{B,(m)}_{n,k}(q) = \sum_{j=0}^{\min(k,n-m-k)} q^{j(j+m)} \qbinom{n+m}{2m+2j}B_{m+j,m}(q) \qbinom{n - 2j-m}{k-j}.
\end{equation}
\end{state}

\begin{proof}
For $0\le k\le\fl{n-m}$, by~\eqref{eq:AB-Nara} and~\eqref{eq:A-positive-expansion},
\[
N^{B,(m)}_{n,k}(q)=\qbinom{n+m}{2m}\qbinom{2m}{m}^{-1}\sum_{j=0}^{\min(k,n-m-k)}q^{j(j+m)}\qbinom{n-m}{2j}A_{m,j}(q)\qbinom{n-2j-m}{k-j}.
\]
For every $j$,
\[
\qbinom{n+m}{2m}\qbinom{n-m}{2j}=\qbinom{n+m}{2m+2j}\qbinom{2m+2j}{2m},
\]
and, directly from the definitions of $A_{m,j}(q)$ in~\eqref{def:sC} and $B_{m+j,m}(q)$ in~\eqref{eq:B-def},
\[
\qbinom{2m+2j}{2m}\qbinom{2m}{m}^{-1}A_{m,j}(q)=B_{m+j,m}(q).
\]
Therefore
\[
\qbinom{n+m}{2m}\qbinom{2m}{m}^{-1}\qbinom{n-m}{2j}A_{m,j}(q)=\qbinom{n+m}{2m+2j}B_{m+j,m}(q),
\]
and substitution into the preceding expansion proves~\eqref{eq:B-positive-expansion} for $0\le k\le\fl{n-m}$. For $\fl{n-m}<k\le n-m$, replace $k$ by $n-m-k$. By~\eqref{eq:B-reflection},
\[
N^{B,(m)}_{n,k}(q)=N^{B,(m)}_{n,n-m-k}(q),
\]
and on the right-hand side of~\eqref{eq:B-positive-expansion} the upper limit is unchanged because
\[
\min(k,n-m-k)=\min(n-m-k,k),
\]
while
\[
\qbinom{n-m-2j}{k-j}=\qbinom{n-m-2j}{n-m-k-j}.
\]
Hence~\eqref{eq:B-positive-expansion} holds for all $0\le k\le n-m$.
\end{proof}

\begin{cor}
\label{state:B-positive}
Let $n$, $m$, and $k$ be nonnegative integers such that $m\leq n$ and $k\leq n-m$. The dual $q$-Narayana numbers $N^{B,(m)}_{n,k}(q)$ are polynomials in $q$ with nonnegative
integer coefficients.
\end{cor}

\begin{proof}
By~\eqref{eq:B-positive-expansion}, $N^{B,(m)}_{n,k}(q)$ is a sum of products of powers of $q$, $q$-binomial coefficients and dual $q$-super Catalan numbers, which have
nonnegative integer coefficients.
\end{proof}

Before turning to $q$-$\gamma$-positivity, we record a second symmetry of the dual refinement: in addition to the row symmetry~\eqref{eq:B-reflection}, each $N^{B,(m)}_{n,k}(q)$ is palindromic as a polynomial in $q$.
\begin{state}
\label{prop:reciprocity}
Let $n$, $m$, and $k$ be nonnegative integers such that $m\leq n$ and $k\leq n-m$, and put $D_{n,m;k}=2\bigl(m(n-m)+k(n-m-k)\bigr)$. Then
\begin{equation}
\label{eq:reciprocity}
q^{D_{n,m;k}}N_{n,k}^{B,(m)}(q^{-1})=N_{n,k}^{B,(m)}(q).
\end{equation}
Consequently, the $N_{n,k}^{B,(m)}(q)$ are monic palindromic polynomials of degree $D_{n,m;k}$.
\end{state}

\begin{proof}
For a $q$-binomial coefficient, $\qbinom{n}{k}\big|_{q\to q^{-1}}=q^{-k(n-k)}\qbinom{n}{k}$. Applied to the four $q$-binomial coefficients in~\eqref{eq:B-Nara}, this gives the
total shift
\[
mn-m^2+k(n-k)+(k+m)(n-m-k)=2\bigl(m(n-m)+k(n-m-k)\bigr).
\]
The numerator and the denominator of~\eqref{eq:B-Nara-factorial} have constant term $1$, so $N_{n,k}^{B,(m)}(q)$ has constant term $1$. By Corollary~\ref{state:B-positive} it
is a polynomial, and reciprocity~\eqref{eq:reciprocity} gives its degree and leading coefficient.
\end{proof}

\section{$q$-$\gamma$-positivity of index one}
\label{sec:qgamma}
Recall Definition~\ref{def:qgamma}: a polynomial $F(q,t)\in\mathbb Z_{\ge0}[q,t]$ of degree $d$ in $t$ is $q$-$\gamma$-positive of index $a$ if
\[
F(q,t)=\sum_{j=0}^{\lfloor d/2\rfloor} q^{a\binom{j+1}{2}}\gamma_j(q)\,t^j \prod_{i=j+1}^{d-j}(1+q^{ai}t)
\qquad\text{with}\qquad \gamma_j(q)\in\mathbb Z_{\ge0}[q].
\]
Now that both $N^{\ast,(m)}_{n,k}(q)$, $\ast\in\{A,B\}$, are known to be polynomials with nonnegative coefficients, we define the weighted $q$-Narayana polynomials corresponding to indices $1$ and $2$ by
\begin{equation}
\label{eq:qrows}
\mathcal{N}^{\ast,(m)}_{n}(q,t)=\sum_{k=0}^{n-m}q^{\binom{k+1}{2}}N^{\ast,(m)}_{n,k}(q)\,t^k, \qquad
\widetilde{\mathcal{N}}^{\ast,(m)}_{n}(q,t)=\sum_{k=0}^{n-m}q^{k(k+1)}N^{\ast,(m)}_{n,k}(q)\,t^k.
\end{equation}
Thus the coefficient of $t^k$ is weighted by $q^{\binom{k+1}{2}}$ for index $1$ and by $q^{k(k+1)}$ for index $2$.

By~\eqref{eq:AB-Nara}, the same duality relation between the two families holds for both $q$-Narayana polynomials:
\begin{equation}
\label{eq:fromAtoB}
\mathcal{N}^{B,(m)}_{n}(q,t)=\frac{B_{n,m}(q)}{A_{n,m}(q)}\,\mathcal{N}^{A,(m)}_{n}(q,t),\qquad
\widetilde{\mathcal{N}}^{B,(m)}_{n}(q,t)=\frac{B_{n,m}(q)}{A_{n,m}(q)}\,\widetilde{\mathcal{N}}^{A,(m)}_{n}(q,t).
\end{equation}
\begin{note}
\label{note:notation}
We compare this notation with the notation in~\cite{RodeletTevlin} and with that of Krattenthaler and Wachs. In~\cite[Eqs.~(75) and~(86)]{RodeletTevlin}, $N^{A,(m)}_{n,k}(q)$ is denoted
$N^{(m)}_{n,k}(q)$, and the classical polynomial $\mathcal{N}^{A,(m)}_n(1,t)$ is denoted $\mathcal{N}^{(m)}_n(1,t)$. The super $q$-Narayana polynomial of Krattenthaler and
Wachs~\cite{Wachs2014} is
\[
N^{(s)}_n(q,t)=\qbinom{2s}{s}\qbinom{n}{s}^{-1}\sum_{r=0}^{n-s}q^{r(r+1)}\qbinom{n}{r}\qbinom{n}{r+s}t^r ;
\]
thus $N^{(s)}_n(q,t)=\widetilde{\mathcal{N}}^{A,(s)}_n(q,t)$ of Sec.~\ref{sec:index-two}, with their parameter $s$ corresponding to our $m$.
\end{note}


We begin with the $q$-$m$-Narayana polynomials for $q$-super Catalans.

\begin{state}
\label{prop:A-qgamma}
Let $m$ and $n$ be nonnegative integers with $m\leq n$. The polynomial $\mathcal{N}^{A,(m)}_{n}(q,t)$ is $q$-$\gamma$-positive of index $1$ and satisfies
\begin{equation}
\label{eq:A-qgamma}
\mathcal{N}^{A,(m)}_{n}(q,t) =\sum_{j=0}^{\fl{n-m}} q^{\binom{j+1}{2}}\gamma^{A,(m)}_{n,j}(q)\, t^j\prod_{i=j+1}^{n-m-j}(1+q^it),
\qquad
\gamma^{A,(m)}_{n,j}(q)=q^{j(j+m)}\qbinom{n-m}{2j}A_{j,m}(q).
\end{equation}
\end{state}

\begin{proof}
By the finite $q$-binomial theorem~\eqref{eq:q-binomial-theorem}, with the parameters $k$ and $n$ there replaced by $j+1$ and $n-m-2j$, respectively,
\[
\prod_{i=j+1}^{n-m-j}(1+q^it) =\sum_{\ell=0}^{n-m-2j} q^{(j+1)\ell+\binom{\ell}{2}}\qbinom{n-m-2j}{\ell}t^\ell .
\]
Comparing coefficients of $t^k$ on both sides of~\eqref{eq:A-qgamma}, with $\ell=k-j$, and using
\begin{equation}
\label{eq:exponent}
\binom{j+1}{2}+(j+1)(k-j)+\binom{k-j}{2}=\binom{k+1}{2},
\end{equation}
the statement is equivalent to
\[
N^{A,(m)}_{n,k}(q)=\sum_{j=0}^{\min(k,n-m-k)} q^{j(j+m)} \qbinom{n-m}{2j} A_{j,m}(q) \qbinom{n-m-2j}{k-j}, \qquad 0\le k\le n-m.
\]
For $k\le\fl{n-m}$ this is exactly~\eqref{eq:A-positive-expansion}, since $A_{m,j}(q)=A_{j,m}(q)$. For larger $k$, both sides are unchanged when $k$ is replaced by $n-m-k$: the left-hand side by~\eqref{eq:B-reflection}; on the right-hand side, the upper limit $\min(k,n-m-k)$ is unchanged and
\[
\qbinom{n-m-2j}{k-j}=\qbinom{n-m-2j}{n-m-k-j}.
\]

The coefficients $\gamma^{A,(m)}_{n,j}(q)$ have nonnegative integer coefficients because each is a product of manifestly nonnegative factors.
\end{proof}

The dual $q$-Narayana polynomials are also $q$-$\gamma$-positive of index $1$.

\begin{state}
\label{prop:B-qgamma}
Let $m$ and $n$ be nonnegative integers with $m\leq n$. The polynomial $\mathcal{N}^{B,(m)}_{n}(q,t)$ is $q$-$\gamma$-positive of index $1$ since it satisfies
\begin{equation}
\label{eq:B-qgamma}
\mathcal{N}^{B,(m)}_{n}(q,t)
=\sum_{j=0}^{\fl{n-m}} q^{\binom{j+1}{2}}\gamma^{B,(m)}_{n,j}(q)\,t^j
\prod_{i=j+1}^{n-m-j}(1+q^it),
\qquad
\gamma^{B,(m)}_{n,j}(q)
=q^{j(j+m)}\qbinom{n+m}{2m+2j}B_{m+j,m}(q).
\end{equation}
\end{state}

\begin{proof}
Comparing coefficients of $t^k$ as in the proof of Proposition~\ref{prop:A-qgamma}, the statement is equivalent to Proposition~\ref{prop:B-positive-expansion}. The coefficients $\gamma^{B,(m)}_{n,j}(q)$ have nonnegative integer coefficients because each is a product of manifestly nonnegative factors.
\end{proof}
\section{$q$-$\gamma$-positivity of index two}
\label{sec:index-two}

For index $2$ we use the polynomials $\widetilde{\mathcal{N}}^{A,(m)}_n(q,t)$ and $\widetilde{\mathcal{N}}^{B,(m)}_n(q,t)$ of~\eqref{eq:qrows}. Krattenthaler and Wachs~\cite{Wachs2014} found the following expression for the type $A$ super $q$-Narayana polynomial; Dilks~\cite{Dilks2015} also treated the case $m=1$:
\begin{equation}
\label{eq:A-index2-gamma}
\widetilde{\mathcal{N}}^{A,(m)}_{n}(q,t) = \sum_{j=0}^{\fl{n-m}}q^{j(j+1)}\,\widetilde{\gamma}^{A,(m)}_{n,j}(q)\,t^j \prod_{i=j+1}^{n-m-j}(1+q^{2i}t),
\end{equation}
where
\[
\widetilde{\gamma}^{A,(m)}_{n,j}(q) = q^j \qbinom{n-m}{2j} A_{j,m} (q^2) \prod_{i=1}^m \frac{1+q^i}{1+q^{m+i}} \prod_{i=1}^{2j} \frac{1 +q^{m+i-1}}{1+q^i}.
\]
They also proved that $\widetilde{\gamma}^{A,(m)}_{n,j}(q)$ is $q$-positive for $m=1$.
However, the coefficients $\widetilde{\gamma}^{A,(m)}_{n,j}(q)$ are not $q$-positive in general. The first example occurs at $m=2$:
\begin{equation}
\label{eq:A-index2-counterexample}
\widetilde{\gamma}^{A,(2)}_{6,2}(q) = q^2\bigl(1+2q^2-q^3+2q^4-q^5+2q^6+q^8\bigr).
\end{equation}
For the dual family, the coefficients $\widetilde{\gamma}^{B,(m)}_{n,j}(q)$ are $q$-positive, and we have the following result:
\begin{state}
\label{prop:B-index2-qgamma}
Let $m$ and $n$ be nonnegative integers with $m\leq n$. The polynomial $\widetilde{\mathcal{N}}^{B,(m)}_{n}(q,t)$ is $q$-$\gamma$-positive of index $2$ since it satisfies
\begin{equation}
\label{eq:B-index2-gamma}
\widetilde{\mathcal{N}}^{B,(m)}_{n}(q,t) =\sum_{j=0}^{\fl{n-m}}q^{j(j+1)}\,\widetilde{\gamma}^{B,(m)}_{n,j}(q)\,t^j \prod_{i=j+1}^{n-m-j}(1+q^{2i}t),
\end{equation}
where
\[
\widetilde{\gamma}^{B,(m)}_{n,j}(q)=q^{j}\qbinom{n+m}{2m+2j} B_{m+j,m}(q^2)\frac{(-q^m;q)_{2j}}{(-q^{2m+1};q)_{2j}}
\]
is a polynomial in $q$ with nonnegative integer coefficients.
\end{state}
This identity follows directly from~\eqref{eq:A-index2-gamma} and~\eqref{eq:fromAtoB}, since
\begin{equation}
\widetilde{\gamma}^{B,(m)}_{n,j}(q) = \qbinom{n+m}{2m} \qbinom{2m}{m}^{-1} \widetilde{\gamma}^{A,(m)}_{n,j}(q).
\end{equation}
The Gaussian factorial formula gives
\[
\qbinom{n+m}{2m}\qbinom{n-m}{2j} \qbinom{2m}{m}^{-1} =\qbinom{n+m}{2m+2j} \qbinom{2m+2j}{2m}\qbinom{2m}{m}^{-1},
\]
while the definitions~\eqref{def:sC} and~\eqref{eq:B-def} give
\[
A_{j,m}(q^2) =\qbinomm{2m}{m}\qbinomm{2m+2j}{2m}^{-1} B_{m+j,m}(q^2).
\]
For Gaussian coefficients, changing the base from $q$ to $q^2$ gives
\[
\qbinomm{2m}{m}\qbinom{2m}{m}^{-1} =\prod_{i=1}^{m}\frac{1+q^{m+i}}{1+q^i}, \qquad \qbinom{2m+2j}{2m}\qbinomm{2m+2j}{2m}^{-1}
=\prod_{i=1}^{2j}\frac{1+q^i}{1+q^{2m+i}}.
\]
Using the expression for $\widetilde{\gamma}^{A,(m)}_{n,j}(q)$ immediately after~\eqref{eq:A-index2-gamma}, the factors indexed by
$1\leq i\leq m$ cancel and the factors indexed by $1\leq i\leq 2j$ simplify as follows:
\begin{align*}
\widetilde{\gamma}^{B,(m)}_{n,j}(q)
&=q^j\qbinom{n+m}{2m+2j}B_{m+j,m}(q^2) \prod_{i=1}^{2j}\frac{1+q^i}{1+q^{2m+i}} \prod_{i=1}^{2j}\frac{1+q^{m+i-1}}{1+q^i}\\
&=q^j\qbinom{n+m}{2m+2j}B_{m+j,m}(q^2)
\frac{(-q^m;q)_{2j}}{(-q^{2m+1};q)_{2j}}.
\end{align*}
Thus the same duality factor continues to relate the two families at index $2$, but it no longer preserves $q$-$\gamma$-positivity: the coefficients on the $A$-side may have negative coefficients, whereas we prove below that the corresponding coefficients on the $B$-side are always $q$-positive.

For completeness, we give an explicit derivation of the expansion~\eqref{eq:B-index2-gamma} in Lemma~\ref{lem:B-index2-identity}. We then prove the positivity of $\widetilde{\gamma}^{B,(m)}_{n,j}(q)$ (Lemma~\ref{lem:X-positive} followed by Lemma~\ref{lem:gamma-positive}).

\subsection{The expansion}
A direct $q$-hypergeometric derivation of the formula for $\widetilde{\gamma}^{B,(m)}_{n,j}(q)$ is based on Singh's quadratic transformation and the $q$-Pfaff--Saalsch\"utz summation. The derivation follows the same steps as those for~\eqref{eq:A-index2-gamma}, which were kindly shared with us by C.~Krattenthaler~\cite{KrattenthalerPC}.

\begin{lem}
\label{lem:B-index2-identity}
Let $m$ and $n$ be nonnegative integers with $m\leq n$. Then the expansion~\eqref{eq:B-index2-gamma} holds.
\end{lem}

\begin{proof}
We prove~\eqref{eq:B-index2-gamma} coefficient by coefficient, in five steps.

\emph{Step 1: reduction to a terminating sum.} The finite $q$-binomial theorem~\eqref{eq:q-binomial-theorem}, with base $q^2$, gives
\[
\prod_{i=j+1}^{n-m-j}(1+q^{2i}t)=\sum_{\ell=0}^{n-m-2j}q^{2(j+1)\ell+\ell(\ell-1)}\qbinomm{n-m-2j}{\ell}t^\ell .
\]
We take $\ell=k-j$. Since
\[
j(j+1)+j+2(j+1)(k-j)+(k-j)(k-j-1)=k(k+1)+j,
\]
the coefficient of $t^k$ on the right-hand side of~\eqref{eq:B-index2-gamma} is $q^{k(k+1)}\sum_{j}c_j$, where
\begin{equation}
\label{eq:cj}
c_j=q^j\qbinom{n+m}{2m+2j}B_{m+j,m}(q^2)\frac{(-q^m;q)_{2j}}{(-q^{2m+1};q)_{2j}}\qbinomm{n-m-2j}{k-j}.
\end{equation}
Therefore~\eqref{eq:B-index2-gamma} is equivalent to
\begin{equation}
\label{eq:index2-coefficients}
N^{B,(m)}_{n,k}(q)=c_0 \sum_{j=0}^{\min(k,n-m-k)}\frac{c_j}{c_0}  \qquad (0\le k\le n-m).
\end{equation}
\emph{Step 2: identification of a terminating ${}_4\phi_3$ series.} 
The initial term in~\eqref{eq:cj} is 
\begin{equation}
\label{eq:c0}
c_0=\qbinom{n+m}{2m}\qbinomm{n-m}{k}
\end{equation}
We transform $c_j/c_0$ and deal with each factor in turn.

Consider the first ratio of binomial factors in $c_j/c_0$ to begin with:
\[
\qbinom{n+m}{2m+2j}\qbinom{n+m}{2m}^{-1}= \frac{(q^{n-m-2j+1};q)_{2m+2j} }{(q;q)_{2m+2j}} \frac{ (q;q)_{2m}}{(q^{n-m+1};q)_{2m}} =  \frac{(q^{n-m-2j+1};q)_{2j}}{(q^{2m+1};q)_{2j}}.
\]
We then combine this ratio with $B_{m+j,m}(q^2)$. By~\eqref{eq:B-def},
\begin{align*}
 B_{m+j,m}(q^2) &=\frac{(q^2;q^2)_{2m+2j}\,(q^2;q^2)_{m}} {(q^2;q^2)_{m+j}\,(q^2;q^2)_{2m}\,(q^2;q^2)_{j}} = \\
&=\frac{(q^2;q^2)_{2m}\,(q^{4m+2};q^2)_{2j}(q^2;q^2)_{m}} {(q^2;q^2)_{m}\,(q^{2m+2};q^2)_{j}(q^2;q^2)_{2m}\,(q^2;q^2)_{j}}= \frac{(q^{4m+2};q^2)_{2j}}{(q^2;q^2)_j\,(q^{2m+2};q^2)_j}.
\end{align*}
Putting these together, we obtain directly from~\eqref{eq:cj}
\begin{align}
\label{eq:step2a}
& \qbinom{n+m}{2m+2j}\qbinom{n+m}{2m}^{-1} B_{m+j,m}(q^2)\frac{(-q^m;q)_{2j}}{(-q^{2m+1};q)_{2j}}=  \notag  \\
&= \frac{(q^{n-m-2j+1};q)_{2j}}{(q^{2m+1};q)_{2j}} \frac{(q^{4m+2};q^2)_{2j}}{(q^2;q^2)_j\,(q^{2m+2};q^2)_j} \frac{(-q^m;q)_{2j}}{(-q^{2m+1};q)_{2j}} = \frac{(-q^m;q)_{2j}(q^{n-m-2j+1};q)_{2j}}{(q^2;q^2)_j(q^{2m+2};q^2)_j} 
\end{align}
where we used \eqref{eq:Ktwo} to combine two terms
\[
(q^{2m+1};q)_{2j}(-q^{2m+1};q)_{2j}=(q^{4m+2};q^2)_{2j}.
\]
Now we simplify the second ratio of binomial factors in $c_j/c_0$. Writing this ratio explicitly and rearranging the factors to make the cancellations visible gives
\begin{align*}
& \qbinomm{n-m-2j}{k-j}\qbinomm{n-m}{k}^{-1} =\frac{(q^2;q^2)_{n-m-2j}}{(q^2;q^2)_{n-m}}  \frac{(q^2;q^2)_{n-m-k}}{(q^2;q^2)_{n-m-k-j}}   \frac{(q^2;q^2)_{k}}{(q^2;q^2)_{k-j}}= \\
& =\frac{(q^{2n-2m-2k-2j+2};q^2)_j (q^{2k-2j+2};q^2)_j }{(q^{2n-2m-4j+2};q^2)_{2j}}
\end{align*}
Going back to \eqref{eq:cj}, we combine this expression with \eqref{eq:step2a} to obtain:
\begin{align*}
& \frac{c_j}{c_0}=q^j \frac{(-q^m;q)_{2j}(q^{n-m-2j+1};q)_{2j}}{(q^2;q^2)_j(q^{2m+2};q^2)_j} \frac{(q^{2n-2m-2k-2j+2};q^2)_j (q^{2k-2j+2};q^2)_j }{(q^{2n-2m-4j+2};q^2)_{2j}} = \\
& = q^j \frac{(-q^m;q)_{2j}(q^{2k-2j+2};q^2)_j(q^{2n-2m-2k-2j+2};q^2)_j}{(q^2;q^2)_j(q^{2m+2};q^2)_j(-q^{n-m-2j+1};q)_{2j}},
\end{align*}
where we have used \eqref{eq:Ktwo} to simplify
\[
(q^{n-m-2j+1};q)_{2j}/(q^{2n-2m-4j+2};q^2)_{2j}=1/(-q^{n-m-2j+1};q)_{2j}.
\]
Now we use the splitting  \eqref{eq:Ksplit}  and reversal  \eqref{eq:Kfour} identities:
\[
(-q^m;q)_{2j}=(-q^m;q^2)_j(-q^{m+1};q^2)_j,
\]
\[
(q^{2k-2j+2};q^2)_j(q^{2n -2m-2k-2j+2};q^2)_j=q^{2(n-m)j-2j(j-1)}(q^{-2k};q^2)_j(q^{-2n + 2m +2k};q^2)_j,
\]
\[
(-q^{n-m-2j+1};q)_{2j}=q^{2(n-m)j-j(2j-1)}(-q^{-n+m};q^2)_j(-q^{1-n+m};q^2)_j.
\]
Consequently,
\begin{equation}
\label{eq:index2-4phi3-terms}
\frac{c_j}{c_0}=q^{2j}\frac{(-q^m;q^2)_j(-q^{m+1};q^2)_j(q^{-2k};q^2)_j(q^{-2n+2m+2k};q^2)_j}{(q^2;q^2)_j(q^{2m+2};q^2)_j(-q^{-n+m};q^2)_j(-q^{1-n+m};q^2)_j}.
\end{equation}
Rewriting~\eqref{eq:index2-coefficients}, we obtain
\begin{equation}
 N^{B,(m)}_{n,k}(q) = c_0 \sum_{j=0}^{\min(k,n-m-k)} q^{2j}\frac{(-q^m;q^2)_j(-q^{m+1};q^2)_j(q^{-2k};q^2)_j(q^{-2n+2m+2k};q^2)_j}{(q^2;q^2)_j(q^{2m+2};q^2)_j(-q^{-n+m};q^2)_j(-q^{1-n+m};q^2)_j}
 \end{equation}
which we identify as a ${}_4\phi_3$ series using~\eqref{eq:series}. Taking $q^2$ as the base and identifying the parameters
\[
 z=q^2, \  (a_1,a_2,a_3,a_4)= \bigl(-q^m,-q^{m+1},q^{-2k},q^{-2n+2m+2k}\bigr), \  (b_1,b_2,b_3) = \bigl(q^{2m+2},-q^{-n+m},-q^{1-n+m}\bigr),
\]
the $j$th summand is exactly the right-hand side of \eqref{eq:index2-4phi3-terms}.  The series is balanced as 
\begin{align*}
&q^2 \,a_1a_2a_3a_4 = q^2 (-q^m)(-q^{m+1})(q^{-2k})(q^{-2n+2m+2k})= q^{4m+3-2n} \quad \text{and} \\
& b_1b_2b_3 = q^{2m+2}(-q^{-n+m})(-q^{1-n+m}) = q^{4m+3-2n}.
\end{align*}
Moreover, the numerator parameters $q^{-2k}$ and $q^{-2n+2m+2k}$ make the series terminating, since
\[
(q^{-2k};q^2)_{k+1}=0 \qquad\text{and}\qquad (q^{-2n+2m+2k};q^2)_{n-m-k+1}=0.
\]
Therefore the series terminates at $ j=\min(k,n-m-k)$.

Thus
\[
 N^{B,(m)}_{n,k}(q) = c_0\, {}_4\phi_3\!\left[
\begin{matrix}
-q^m,-q^{m+1},q^{-2k},q^{-2n +2m+2k}\\
q^{2m+2},-q^{-n+m},-q^{1-n+m}
\end{matrix};
q^2,q^2
\right].
\]
Using the value of $c_0$ from \eqref{eq:c0} gives
\begin{equation}
\label{eq:index2-4phi3}
N^{B,(m)}_{n,k}(q)=\qbinom{n+m}{2m}\qbinomm{n-m}{k}\,{}_4\phi_3\!\left[\begin{matrix}-q^m,-q^{m+1},q^{-2k},q^{-2n +2m+2k}\\ q^{2m+2},-q^{-n+m},-q^{1-n+m}\end{matrix};q^2,q^2\right],
\end{equation}
where the series terminates at $j=\min(k,n-m-k)$.

\emph{Step 3: Singh's quadratic transformation.} The next crucial step is Singh's quadratic transformation~\cite[Appendix~(III.21)]{GR}, which we record here for reference:
\begin{equation}
\label{eq:Singh}
{}_4\phi_3\!\left[\begin{matrix} a^2,b^2,c,d\\ abq^{1/2},-abq^{1/2},-cd\end{matrix};q,q\right] ={}_4\phi_3\!\left[\begin{matrix} a^2,b^2,c^2,d^2\\a^2b^2q,-cd,-cdq\end{matrix};q^2,q^2\right],
\end{equation}
which is valid when the series terminate. The series in~\eqref{eq:index2-4phi3} has the form of the right-hand side of~\eqref{eq:Singh}. Indeed, we choose
\begin{itemize}
\item $a^2=-q^{m+1}$ and $b^2=-q^m$, so that $a^2b^2q=q^{2m+2}$ and $\{abq^{1/2},-abq^{1/2}\}=\{q^{m+1},-q^{m+1}\}$, whichever square root is taken;
\item $c=q^{k-n+m}$ and $d=q^{-k}$, so that $c^2=q^{2k-2n+2m}$, $d^2=q^{-2k}$, $-cd=-q^{-n+m}$ and $-cdq=-q^{1-n+m}$.
\end{itemize}
Thus~\eqref{eq:index2-4phi3} becomes
\begin{equation}
\label{eq:index2-4phi3-singh}
N^{B,(m)}_{n,k}(q)=\qbinom{n+m}{2m}\qbinomm{n-m}{k} {}_4\phi_3\!\left[
\begin{matrix} -q^{m+1},-q^m,q^{k-n+m},q^{-k}\\q^{m+1},-q^{m+1},-q^{-n+m}
\end{matrix};q,q\right].
\end{equation}
The next simplification is quite remarkable. Notice that the numerator parameter $-q^{m+1}$ can be cancelled against the equal denominator parameter, so 
that in~\eqref{eq:index2-4phi3-singh} the ${}_4\phi_3$ series reduces to a ${}_3\phi_2$ series, and
\begin{align}
& N^{B,(m)}_{n,k}(q)=\qbinom{n+m}{2m}\qbinomm{n-m}{k}{}_4\phi_3\!\left[\begin{matrix} -q^{m+1},-q^m,q^{k-n+m},q^{-k}\\ q^{m+1},-q^{m+1},-q^{-n+m}\end{matrix};q,q\right] \notag \\
\label{eq:3phi2}
& = \qbinom{n+m}{2m}\qbinomm{n-m}{k} {}_3\phi_2\!\left[\begin{matrix} -q^m,q^{k-n+m},q^{-k}\\ q^{m+1},-q^{-n+m}\end{matrix};q,q\right].
\end{align}
\emph{Step 4: the $q$-Pfaff--Saalsch\"utz summation.} To make further progress, we will need the $q$-Pfaff--Saalsch\"utz summation~(II.12) of~\cite{GR}, which we rewrite here:
\begin{equation}
\label{eq:Saal}
{}_3\phi_2\!\left[\begin{matrix} a,b,q^{-N}\\c,abq^{1-N}/c\end{matrix};q,q\right]=\frac{(c/a;q)_N\,(c/b;q)_N}{(c;q)_N\,(c/(ab);q)_N}.
\end{equation}
The ${}_3\phi_2$ series obtained in Step~3 has this form, with
\begin{itemize}
\item $N=k$, $a=-q^m$, $b=q^{k-n+m}$ and $c=q^{m+1}$;
\item $abq^{1-N}/c=-q^{m}q^{k-n+m}q^{1-k}/q^{m+1}=-q^{-n+m}$, as required;
\item $c/a=-q$, $c/b=q^{n-k+1}$ and $c/(ab)=-q^{n-m-k+1}$.
\end{itemize}
Hence
\[
{}_3\phi_2\!\left[\begin{matrix} -q^m,q^{k-n+m},q^{-k}\\ q^{m+1},-q^{-n+m}\end{matrix};q,q\right] =\frac{(-q;q)_k\,(q^{n-k+1};q)_k}{(q^{m+1};q)_k\,(-q^{n-m-k+1};q)_k}.
\]
And
\begin{equation}
\label{eq:index2-final}
N^{B,(m)}_{n,k}(q)=\qbinom{n+m}{2m}\qbinomm{n-m}{k} \frac{(-q;q)_k\,(q^{n-k+1};q)_k}{(q^{m+1};q)_k\,(-q^{n-m-k+1};q)_k}.
\end{equation}
\emph{Step 5: final simplification.} For the final step, we use~\eqref{eq:Ktwo} to rewrite
\[
\qbinomm{n-m}{k}=\frac{(q^{2n-2m-2k+2};q^2)_k}{(q^2;q^2)_k}=\frac{(q^{n-m-k+1};q)_k\,(-q^{n-m-k+1};q)_k}{(q;q)_k\,(-q;q)_k}.
\]
Thus~\eqref{eq:index2-final} takes the form

\begin{align*}
& N^{B,(m)}_{n,k}(q) = \qbinom{n+m}{2m}\qbinomm{n-m}{k} \frac{(-q;q)_k\,(q^{n-k+1};q)_k}{(q^{m+1};q)_k\,(-q^{n-m-k+1};q)_k} = \\
&= \frac{(q^{n-m+1};q)_{2m}}{(q;q)_{2m}} \frac{(q^{n-m-k+1};q)_k\,(-q^{n-m-k+1};q)_k}{(q;q)_k\,(-q;q)_k} \frac{(-q;q)_k\,(q^{n-k+1};q)_k}{(q^{m+1};q)_k\,(-q^{n-m-k+1};q)_k} = \\
&= \frac{(q^{n-m+1};q)_{2m}\,(q^{n-m-k+1};q)_k\,(q^{n-k+1};q)_k} {(q;q)_{2m}\,(q;q)_k\,(q^{m+1};q)_k} = \frac{(q;q)_{n+m}\,(q;q)_m\,(q;q)_n} {(q;q)_{2m}\,(q;q)_{n-m-k}\,(q;q)_k\,(q;q)_{n-k}\,(q;q)_{m+k}},
\end{align*}

which is exactly~\eqref{eq:B-Nara-factorial}. This proves~\eqref{eq:index2-coefficients}.
\end{proof}
\subsection{Positivity of the $\gamma$ coefficients}
It remains to show that the rational functions $\widetilde{\gamma}^{B,(m)}_{n,j}(q)$ are polynomials with nonnegative integer coefficients.  
Let's recall that
\[
\widetilde{\gamma}^{B,(m)}_{n,j}(q) = q^{j}\qbinom{n+m}{2m+2j} B_{m+j,m}(q^2)\frac{(-q^m;q)_{2j}}{(-q^{2m+1};q)_{2j}}.
\]
Since $q^{j}\qbinom{n+m}{2m+2j} $ is manifestly $q$-positive, we need to show the second factor is also $q$-positive.
Let's use the notation
\begin{equation}
\label{eq:X-def}
X_{m,j}(q)=B_{m+j,m}(q^2)\frac{(-q^m;q)_{2j}}{(-q^{2m+1};q)_{2j}} = \frac{(-q^m;q)_{2j}\,(q^{2m+1};q^2)_j}{(q^2;q^2)_j}
\qquad(m,j\ge0)
\end{equation}
Our immediate goal is to derive the boundary values and establish a recurrence for $X_{m,j}(q)$ with positive coefficients.

\begin{lem}
\label{lem:X-base}
The boundary values of~\eqref{eq:X-def} are
\begin{equation}
\label{eq:X-initial}
X_{m,0}(q)=1\quad(m\ge0),\qquad
X_{0,j}(q)=2\qbinomm{2j-1}{j}\quad(j\ge1).
\end{equation}
\end{lem}
\begin{proof}
For $j=0$, the products in~\eqref{eq:X-def} are empty, so $X_{m,0}(q)=1$. For $m=0$ and $j\ge1$,
\[
X_{0,j}(q) = \frac{(-1;q)_{2j}(q;q^2)_j}{(q^2;q^2)_j} = \frac{ 2(-q;q^2)_j(-q^2;q^2)_{j-1}(q;q^2)_j }{(q^2;q^2)_j },
\]
since
\[
(-1;q)_{2j} = 2(-q;q)_{2j-1} = 2(-q;q^2)_j(-q^2;q^2)_{j-1}.
\]
Using~\eqref{eq:Ktwo} twice,
\[
(-q;q^2)_j(q;q^2)_j=(q^2;q^4)_j \quad \text{and} \quad  (-q^2;q^2)_{j-1} = \frac{(q^4;q^4)_{j-1}}{(q^2;q^2)_{j-1}},
\]
this becomes
\begin{equation}
X_{0,j}(q) =
\frac{ 2(q^2;q^4)_j(q^4;q^4)_{j-1} }{ (q^2;q^2)_j(q^2;q^2)_{j-1} }= \frac{ 2(q^2;q^2)_{2j-1} }{ (q^2;q^2)_j(q^2;q^2)_{j-1} }= 2\qbinomm{2j-1}{j},
\end{equation}
where in the second equality we used
\[
(q^2;q^2)_{2j-1} =(q^2;q^4)_j(q^4;q^4)_{j-1}.
\]
\end{proof}

We now establish a recurrence for $X_{m,j}(q)$:
\begin{lem}
\label{lem:X-recurrence}
For every $m,j\ge1$, the rational functions $X_{m,j}(q)$ defined in~\eqref{eq:X-def} satisfy
\begin{equation}
\label{eq:X-positive-recurrence}
X_{m,j}(q)=q^{2j}X_{m-1,j}(q)+(1+q^{m+2j-2})(1+q^{3m+2j-2})X_{m,j-1}(q).
\end{equation}
\end{lem}
\begin{proof}
From~\eqref{eq:X-def}, for $m,j\ge1$,
\[
\frac{X_{m,j}(q)}{X_{m,j-1}(q)}
=\frac{(1+q^{m+2j-2})(1+q^{m+2j-1})(1-q^{2m+2j-1})}{1-q^{2j}}
\]
and
\[
\frac{X_{m-1,j}(q)}{X_{m,j-1}(q)}
=\frac{(1+q^{m-1})(1+q^{m+2j-2})(1-q^{2m-1})}{1-q^{2j}}.
\]
Hence
\[
\frac{X_{m,j}(q)-q^{2j}X_{m-1,j}(q)}{X_{m,j-1}(q)}
=\frac{1+q^{m+2j-2}}{1-q^{2j}}
\Bigl((1+q^{m+2j-1})(1-q^{2m+2j-1})-q^{2j}(1+q^{m-1})(1-q^{2m-1})\Bigr).
\]
Expanding, the expression in parentheses is
\begin{align*}
&1+q^{m+2j-1}-q^{2m+2j-1}-q^{3m+4j-2}-q^{2j}-q^{m+2j-1}+q^{2m+2j-1}+q^{3m+2j-2}\\
&=(1-q^{2j})(1+q^{3m+2j-2}),
\end{align*}
which proves~\eqref{eq:X-positive-recurrence}. 
\end{proof}
\begin{lem}[Positivity of $X_{m,j}(q)$]
\label{lem:X-positive}
\begin{equation}
\label{eq:X-pos}
X_{m,j}(q) \in \mathbb Z_{\ge0}[q], \quad m, j \geq 0
\end{equation}
\end{lem}
\begin{proof}
We induct on $m+j$.  The boundary cases $j=0$ and $m=0$ follow from Lemma~\ref{lem:X-base}. Let $m,j\ge1$, and suppose the assertion holds for all pairs
$(m',j')$ with $m'+j'<m+j$.  By~\eqref{eq:X-positive-recurrence},
\[
X_{m,j}(q) = q^{2j}X_{m-1,j}(q) + (1+q^{m+2j-2})(1+q^{3m+2j-2})X_{m,j-1}(q).
\]
The induction hypothesis gives
\[
X_{m-1,j}(q),\;X_{m,j-1}(q)\in\mathbb Z_{\ge0}[q].
\]
The coefficients multiplying these two polynomials in~\eqref{eq:X-positive-recurrence} have nonnegative coefficients. Hence both terms on the right-hand side belong to $\mathbb Z_{\ge0}[q]$, and therefore so does $X_{m,j}(q)$.
\end{proof}
\begin{lem}[Positivity of $\widetilde{\gamma}^{B,(m)}_{n,j}(q)$]
\label{lem:gamma-positive}
\begin{equation}
\label{eq:gamma-pos}
\widetilde{\gamma}^{B,(m)}_{n,j}(q) \in \mathbb Z_{\ge0}[q] \quad  0 \le m \leq n, \ 0 \leq j \leq \fl{n-m}
\end{equation}
\end{lem}
\begin{proof}
By Lemma~\ref{lem:X-positive}, $X_{m,j}(q)\in\mathbb Z_{\ge0}[q]$. Hence $\widetilde{\gamma}^{B,(m)}_{n,j}(q)=q^{j}\qbinom{n+m}{2m+2j}X_{m,j}(q)$ is a polynomial with nonnegative integer coefficients.
\end{proof}
Collecting all the results in this section, we have
\begin{proof}[Proof of Proposition~\ref{prop:B-index2-qgamma}]
The expansion~\eqref{eq:B-index2-gamma} is Lemma~\ref{lem:B-index2-identity}, as required by Definition~\ref{def:qgamma} with $a=2$. The coefficients $\widetilde{\gamma}^{B,(m)}_{n,j}(q)$ have nonnegative integer coefficients by Lemma~\ref{lem:gamma-positive}.
\end{proof}
\section{$q$-Touchard Identities}
\label{sec:touchard}
In~\cite[Proposition~1]{RodeletTevlin} we extended Andrews' $q$-Touchard identity to all $q$-super Catalan numbers:
\[
A_{n+m,m}(q) = \sum_{k = 0}^{\fl{n}} q^{2k(k+m)} \qbinom{n}{2k} \frac{(-q^{k+ m+1};q)_{n-k}}{(-q;q)_k} A_{k,m}(q).
\]
Using the duality relation~\eqref{eq:AB-ratio}, together with the symmetry $A_{k,m}(q)=A_{m,k}(q)$, the preceding identity becomes:
\begin{equation}
\label{eq:B-Touchard}
B_{n+m,m}(q) = \sum_{k = 0}^{\fl{n}} q^{2k(k+m)} \qbinom{n+2m}{2k+2m} \frac{(-q^{k+m+1};q)_{n-k}}{(-q;q)_k} B_{k+m,m}(q).
\end{equation}
Combining~\eqref{eq:B-weighted} with~\eqref{eq:B-positive-expansion}, we get another $q$-Touchard identity:
\begin{state}[$q$-Touchard identity]
\label{prop:B-Touchard-second}
For all nonnegative integers $m$ and $n$,
\begin{equation}
\label{eq:B-Touchard-second}
B_{n+m,m}(q)=\sum_{j=0}^{\fl n}q^{2j(j+m)}\qbinom{n+2m}{2m+2j}B_{m+j,m}(q)\sum_{\ell=0}^{n-2j}q^{\ell(\ell+2j+m)}\qbinom{n-2j}{\ell}.
\end{equation}
At  $q=1$ both forms specialize to the classical expression
\begin{equation}
\label{eq:B-Touchard-classical}
B_{n+m,m}=\sum_{j=0}^{\fl n}2^{\,n-2j}\binom{n+2m}{2m+2j}B_{m+j,m}.
\end{equation}
\end{state}
The first identity~\eqref{eq:B-Touchard} is transported from the $A$-family by duality, whereas~\eqref{eq:B-Touchard-second} arises directly from the dual $q$-Narayana refinement.
\begin{proof}
By~\eqref{eq:B-weighted},
\[
B_{n+m,m}(q)=\sum_{k=0}^{n}q^{k(k+m)}N^{B,(m)}_{n+m,k}(q).
\]
Applying~\eqref{eq:B-positive-expansion} with $n+m$ in place of $n$ gives
\[
N^{B,(m)}_{n+m,k}(q)=\sum_{j=0}^{\min(k,n-k)}q^{j(j+m)}\qbinom{n+2m}{2m+2j}B_{m+j,m}(q)\qbinom{n-2j}{k-j}.
\]
Substitution therefore yields
\[
B_{n+m,m}(q)=\sum_{k=0}^{n}\sum_{j=0}^{\min(k,n-k)}q^{k(k+m)+j(j+m)}\qbinom{n+2m}{2m+2j}B_{m+j,m}(q)\qbinom{n-2j}{k-j}.
\]
Since $0\le j\le\min(k,n-k)$ is equivalent to $0\le j\le\fl n$ and $j\le k\le n-j$, reversing the order of summation gives
\[
B_{n+m,m}(q)=\sum_{j=0}^{\fl n}q^{j(j+m)}\qbinom{n+2m}{2m+2j}B_{m+j,m}(q)\sum_{k=j}^{n-j}q^{k(k+m)}\qbinom{n-2j}{k-j}.
\]
Set $\ell=k-j$. Then $0\le\ell\le n-2j$ and
\[
k(k+m)+j(j+m)=(\ell+j)(\ell+j+m)+j(j+m)=\ell(\ell+2j+m)+2j(j+m).
\]
Consequently,
\[
B_{n+m,m}(q)=\sum_{j=0}^{\fl n}q^{2j(j+m)}\qbinom{n+2m}{2m+2j}B_{m+j,m}(q)\sum_{\ell=0}^{n-2j}q^{\ell(\ell+2j+m)}\qbinom{n-2j}{\ell},
\]
which proves~\eqref{eq:B-Touchard-second}
\end{proof}
\section{$q$-unimodality of $q$-Narayana numbers}
\label{sec:unimodality}

For polynomials $f(q)$ and $g(q)$, write $f(q)\le_q g(q)$ if $g(q)-f(q)\in\mathbb Z_{\ge0}[q]$. Krattenthaler and Wachs~\cite{Wachs2014} proved the coefficient-wise $q$-unimodality of the rows $N^{A,(m)}_{n,k}(q)$; here we recover this result and prove the corresponding statement for $N^{B,(m)}_{n,k}(q)$.

The inequalities in Corollary~\ref{cor:row-q-unimodal} can also be proved directly from Propositions~\ref{prop:A-qgamma} and~\ref{prop:B-qgamma} by comparing consecutive coefficients and using the coefficient-wise positivity of the differences between consecutive Gaussian coefficients. 

We choose instead to deduce the result from the following consequence of $q$-$\gamma$-positivity due to Krattenthaler and Wachs~\cite{Wachs2014}, since this gives the symmetry and coefficient-wise $q$-unimodality of both families directly from their $q$-$\gamma$-positivity of index $1$.
 If
\[
P(q,t)=\sum_{k=0}^{d}p_k(q)t^k\in\mathbb Z_{\ge0}[q,t]
\]
is $q$-$\gamma$-positive of index $a$, then the polynomial obtained by multiplying the coefficient of $t^k$ by $q^{-a \binom{k+1}{2}}$,
\[
\widetilde P(q,t)=\sum_{k=0}^{d}q^{-a\binom{k+1}{2}}p_k(q)t^k,
\]
is palindromic and coefficient-wise $q$-unimodal.

\begin{cor}
\label{cor:row-q-unimodal}
Let $m$ and $n$ be nonnegative integers with $m\le n$. Then the rows of $q$-$m$-Narayana numbers and dual $q$-Narayana numbers are symmetric and coefficient-wise $q$-unimodal. Thus, for $0\le 2k<n-m$,
\begin{equation}
\label{eq:row-difference}
N^{A,(m)}_{n,k+1}(q)-N^{A,(m)}_{n,k}(q)\in\mathbb Z_{\ge0}[q],
\qquad
N^{B,(m)}_{n,k+1}(q)-N^{B,(m)}_{n,k}(q)\in\mathbb Z_{\ge0}[q],
\end{equation}
while for $0\le k\le n-m$,
\begin{equation}
\label{eq:row-symmetry}
N^{A,(m)}_{n,k}(q)=N^{A,(m)}_{n,n-m-k}(q),
\qquad
N^{B,(m)}_{n,k}(q)=N^{B,(m)}_{n,n-m-k}(q).
\end{equation}
\end{cor}

\begin{proof}
By Proposition~\ref{prop:A-qgamma},
\[
\mathcal N^{A,(m)}_n(q,t)=\sum_{k=0}^{n-m}q^{\binom{k+1}{2}}N^{A,(m)}_{n,k}(q)t^k
\]
is $q$-$\gamma$-positive of index $1$. Applying the Krattenthaler--Wachs theorem with $a=1$, the coefficient of $t^k$ is multiplied by $q^{-\binom{k+1}{2}}$, and hence
\[
\sum_{k=0}^{n-m}N^{A,(m)}_{n,k}(q)t^k
\]
is palindromic and coefficient-wise $q$-unimodal.

Similarly, Proposition~\ref{prop:B-qgamma} gives that
\[
\mathcal N^{B,(m)}_n(q,t)=\sum_{k=0}^{n-m}q^{\binom{k+1}{2}}N^{B,(m)}_{n,k}(q)t^k
\]
is $q$-$\gamma$-positive of index $1$. Again multiplying the coefficient of $t^k$ by $q^{-\binom{k+1}{2}}$ gives
\[
\sum_{k=0}^{n-m}N^{B,(m)}_{n,k}(q)t^k
\]
is also palindromic and coefficient-wise $q$-unimodal. The inequalities~\eqref{eq:row-difference} and symmetries~\eqref{eq:row-symmetry} are exactly these two conclusions written coefficient-wise.
\end{proof}

At $q=1$, we have

\begin{cor}
\label{cor:B-classical-gamma}
Let $m$ and $n$ be nonnegative integers with $m\leq n$. The dual Narayana polynomial $\mathcal{N}^{B,(m)}_{n}(1,t)$ is $\gamma$-positive since it satisfies
\begin{equation}
\label{eq:B-classical-gamma}
\mathcal{N}^{B,(m)}_{n}(1,t) =\sum_{j=0}^{\fl{n-m}} \binom{n+m}{2m+2j}B_{m+j,m} \,t^j(1+t)^{n-m-2j}.
\end{equation}
In particular, the dual Narayana numbers $N^{B,(m)}_{n,k}$, $0\le k\le n-m$, form a symmetric unimodal sequence.
\end{cor}

\begin{proof}
Put $q=1$ in~\eqref{eq:B-qgamma}. The coefficients are nonnegative integers and every polynomial $t^j(1+t)^{n-m-2j}$ has a symmetric unimodal coefficient sequence centered at $(n-m)/2$.
\end{proof}

The corresponding identity for the $m$-Narayana polynomials, $\mathcal{N}^{A,(m)}_{n}(1,t)=\sum_j\binom{n-m}{2j}A_{j,m}\,t^j(1+t)^{n-m-2j}$, is~\cite[Eq.~(76)]{RodeletTevlin}.
Note that the binomial coefficient there has top index $n-m$, whereas in~\eqref{eq:B-classical-gamma} it has top index $n+m$.

\section{Acknowledgements}
The author would like to thank Christian Krattenthaler for sharing his notes and Jean-Christophe Novelli for many useful discussions.

\bibliographystyle{amsplain}
\bibliography{dual_super_catalan}

\end{document}